# A posteriori error analysis for a boundary element method with non-conforming domain decomposition [*]

Catalina Domínguez [†]      Norbert Heuer [‡]

## Abstract

We present and analyze an a posteriori error estimator based on mesh refinement for the solution of the hypersingular boundary integral equation governing the Laplacian in three dimensions. The discretization under consideration is a non-conforming domain decomposition method based on the Nitsche technique. Assuming a saturation property, we establish quasi-reliability and efficiency of the error estimator in comparison with the error in a natural (non-conforming) norm. Numerical experiments with uniform and adaptively refined meshes confirm our theoretical results.

## 1  Introduction

In the last years we have started to design and analyze non-conforming Galerkin approximations for solutions to hypersingular boundary integral equations in three dimensions. There are non-conforming sub-domain based variants like the mortar method [15] and a Nitsche coupling [8]. Element-wise discontinuous approximations are considered in [22] (Crouzeix-Raviart elements) and [20] (Nitsche-based discontinuous Galerkin method). Main advantage of these methods is their flexibility in the construction of discrete spaces, e.g., on complicated surfaces or for ease of adaptivity. In this paper we present and study an a posteriori error estimator for non-conforming boundary elements, more precisely for the Nitsche domain decomposition method presented in [8].

The so-called energy space for the variational setting of hypersingular integral equations is a Sobolev space of order $1/2$. Principal difficulty in the analysis of their non-conforming (discontinuous in this case) discretization is that there is no well-defined trace operator in the energy space. Such an operator is needed to deal with jump conditions in variational form. Therefore, a priori analyses of these methods avoid the use of variational settings at the continuous level.

---

[*]Supported by CONICYT through FONDECYT project 1110324 and Anillo ACT1118 (ANANUM), and by Agenda IDI-2012 Uninorte.

[†]Departamento de Matemáticas y Estadística, Universidad del Norte, Barranquilla, Colombia, email: `dcatalina@uninorte.edu.co`

[‡]Facultad de Matemáticas, Pontificia Universidad Católica de Chile, Avenida Vicuña Mackenna 4860, Macul, Santiago, Chile, email: `nheuer@mat.puc.cl`

Instead, a technique is used which gives rise to the second Strang lemma, i.e., discrete ellipticity and continuity (for more regular functions) of the principal bilinear form. In the case of element-wise discontinuous functions, main ingredients are scaling properties of fractional-order Sobolev (semi-) norms as well as arguments from the equivalence of norms in finite-dimensional spaces. Such arguments are not available when studying approximations which are discontinuous across a sub-domain decomposition (since discrete spaces on individual sub-domains do not have bounded dimension when meshes are refined). Therefore, the analysis of non-conforming sub-domain based methods ignores scaling arguments (sub-domains are considered fixed) and rather uses techniques from Sobolev spaces (of higher regularity than $1/2$) and the inverse property of discrete functions. By these arguments, (poly-) logarithmic perturbations appear in a priori error estimates, cf. [15, 8].

A posteriori error analysis faces the same problems, and requires different approaches for element and sub-domain-wise discontinuous approximations. In [19] we study an error estimator based on mesh refinement for Crouzeix-Raviart boundary elements. The result is somewhat surprising in the sense that adaptivity is advantageous (convergence is faster) only when solutions are more singular than usual, i.e., when there are stronger than usual (edge) singularities. In this paper we study a posteriori error estimation for the non-conforming boundary element method based on Nitsche domain decomposition [8]. Contrary to adaptive Crouzeix-Raviart boundary elements, here we observe improved convergence orders when using adaptivity for model problems with standard edge singularities. Due to the difficulties previously discussed, efficiency can be proved but reliability is slightly perturbed by an $\epsilon$-reduced order of convergence (with any $\epsilon > 0$) or by a factor that blows up logarithmically when the mesh size goes to zero. As mentioned, such perturbations appear also in a priori error estimates. Nevertheless, these effects are difficult to observe in numerical experiments.

Let us mention some connections with a posteriori error estimates for conforming boundary elements. Adaptive methods based on mesh refinement (halving mesh sizes) are also called $h$-$h/2$ methods. Their reliability usually depends on a saturation assumption, cf. [13, 10]. Therefore, it is not surprising that we need this assumption here as well. Note, however, that averaging techniques also suffer from such a restriction (reliability can be shown up to terms which are considered to be of higher order), cf. [7] for the hypersingular operator for problems in two dimensions. Very similar to $h$-$h/2$ methods are strategies based on space enrichments combined with stable additive Schwarz decompositions, cf. [25, 21, 17, 11]. Here, too, a saturation assumption is needed to prove reliability. The situation is different with residual based estimators which do without such an assumption, cf., e.g., [5, 6]. However, it is an open problem how to analyze residual based estimators for non-conforming boundary elements. Principal reason for this is that it is not clear how to define a residual for the hypersingular integral equation when discontinuous approximations are used since, in this case, the operator is not well defined.

The previously mentioned papers on adaptive conforming boundary elements all deal with symmetric positive definite operators. In this paper, with a domain decomposition setting, the underlying bilinear form is not symmetric so that orthogonal projection arguments do not apply. Therefore, the saturation assumption also enters the efficiency estimate, though in a weaker form than it influences reliability. This is known from finite element a posteriori error analysis of indefinite operators, cf. [1], and also from two-level boundary element estimators for acoustic scattering (which lead to indefinite operators), cf. [23].

In this paper we use techniques from a priori error analysis of non-conforming boundary elements [8], from a posteriori error analysis of conforming boundary elements [13], estimates for a Clement type interpolation operator [4], and inverse properties of piecewise polynomials in fractional-order Sobolev spaces [16]. The main result (Theorem 4.1) establishes efficiency and quasi-reliability in a natural norm. This natural norm (denoted $\|\cdot\|_a$ or $a$-norm below) does not decompose the appearing surface differential operator (**curl**) within its Sobolev space of order $-1/2$ and, thus, avoids some complications. In a corollary (Corollary 4.6) we illustrate what can be shown in a more standard (broken) Sobolev norm. Here, quasi-reliability is established in a broken Sobolev norm of order $1/2$ (this is the expected order for a hypersingular operator) but efficiency is only obtained for the error measured in a norm of slightly higher order.

The remainder of the paper is organized as follows. In Section 2 we define some norms and present the model problem. The non-conforming boundary element method with Nitsche coupling is recalled in Section 3. In Section 4 we present the error estimator, the main result (Theorem 4.1), several technical details (Subsection 4.1), a proof of Theorem 4.1 (Subsection 4.2), the error estimate in standard broken norms (Corollary 4.6 in Subsection 4.3), and a standard adaptive refinement algorithm (Subsection 4.4) to be used in our numerical experiments. These experiments are presented in Section 5.

## 2 Preliminaries and model problem

First, we recall the definition of some fractional-order Sobolev norms, see, e.g., [24]. For a Lipschitz domain $\Omega \subset \mathbb{R}^n$ and $0 < s < 1$ we define

$$\|u\|^2_{H^s(\Omega)} := \|u\|^2_{L^2(\Omega)} + |u|^2_{H^s(\Omega)}$$

where

$$|u|_{H^s(\Omega)} := \left( \int_\Omega \int_\Omega \frac{|u(x) - u(y)|^2}{|x - y|^{2s+n}} \, dx \, dy \right)^{1/2}$$

is a semi-norm in $H^s(\Omega)$. For $0 < s < 1$, the space $\tilde{H}^s(\Omega)$ is defined as the completion of $C_0^\infty(\Omega)$ under the norm

$$\|u\|_{\tilde{H}^s(\Omega)} := \left( |u|^2_{H^s(\Omega)} + \int_\Omega \frac{u(x)^2}{(\mathrm{dist}(x, \partial\Omega))^{2s}} \, dx \right)^{1/2}.$$

When $s \in (0, 1/2)$, $\|\cdot\|_{\tilde{H}^s(\Omega)}$ and $\|\cdot\|_{H^s(\Omega)}$ are equivalent norms whereas for $s \in (1/2, 1)$ there holds $\tilde{H}^s(\Omega) = H_0^s(\Omega)$, the latter space consisting of $H^s(\Omega)$-function with vanishing trace on the boundary of $\Omega$. For $s > 0$ the spaces $H^{-s}(\Omega)$ and $\tilde{H}^{-s}(\Omega)$ are the dual spaces (with $L^2(\Omega)$ as pivot space) of $\tilde{H}^s(\Omega)$ and $H^s(\Omega)$, respectively.

Throughout the paper, the notation $\langle \cdot, \cdot \rangle_\Gamma$ denotes the duality pairing extending the $L^2(\Gamma)$-inner product. Also, $a \lesssim b$ means that $a \leq cb$ with a generic constant $c > 0$ that is independent of involved parameters like $h$. Similarly, the notation $a \gtrsim b$ and $a \simeq b$ is used.

Now let $\Gamma$ be a plane open surface with polygonal boundary $\partial \Gamma$. We will identify sub-domains of $\Gamma$ with polygonal subsets of $\mathbb{R}^2$. Our model problem is as follows. *For a given function $f \in L^2(\Gamma)$ find $u \in \tilde{H}^{1/2}(\Gamma)$ such that*

$$Wu := -\frac{1}{4\pi} \frac{\partial}{\partial \mathbf{n}} \int_\Gamma u(y) \frac{\partial}{\partial \mathbf{n}_y} \frac{1}{|\cdot -y|} dS_y = f \quad on \quad \Gamma. \tag{1}$$

Here, $W$ is the hypersingular integral operator of the Laplacian in three space dimensions and $\mathbf{n}$ is the exterior normal unit vector on $\Gamma$ in a certain direction. Equation (1) represents a Neumann problem for the Laplacian in $\mathbb{R}^3 \setminus \bar{\Gamma}$.

A weak form of (1) is: *Find $u \in \tilde{H}^{1/2}(\Gamma)$ such that*

$$\langle Wu, v \rangle_\Gamma = \langle f, v \rangle_\Gamma \quad \forall v \in \tilde{H}^{1/2}(\Gamma). \tag{2}$$

This is a coercive problem since $W : \tilde{H}^{1/2}(\Gamma) \to H^{-1/2}(\Gamma)$ is continuous and coercive, cf. [9].
A conforming boundary element method (BEM) for the approximate solution of (2) is

$$\tilde{u}_h \in \tilde{H}_h : \quad \langle W\tilde{u}_h, v \rangle_\Gamma = \langle f, v \rangle_\Gamma \quad \forall v \in \tilde{H}_h \tag{3}$$

with $\tilde{H}_h \subset \tilde{H}^{1/2}(\Gamma)$ being a piecewise polynomial subspace. This conformity condition requires that approximating functions be continuous with zero trace on $\partial \Gamma$.

In the following we consider the non-conforming domain decomposition method from [8] and analyze an a posteriori error estimation.

## 3 Boundary element method with non-conforming domain decomposition

In this section we recall the discrete domain decomposition setting from [8] for a non-conforming solution of our model problem (1).

Let $\mathcal{T}$ be a decomposition of $\Gamma$ into non-overlapping polygonal sub-domains (sub-surfaces) $\Gamma_1, \ldots, \Gamma_K$. The (relatively closed) interface between two neighboring sub-domains $\Gamma_i$, $\Gamma_j$ with $i \neq j$ is denoted by $\gamma_{ij}$. We assume that $\gamma_{ij}$ is an entire edge of $\Gamma_i$ or $\Gamma_j$ (other situations like an interface reducing to a single point or subsets of edges are not allowed). The skeleton of $\mathcal{T}$ is

$$\gamma := \cup \{\gamma_{ij}; \ \bar{\Gamma}_i \cap \bar{\Gamma}_j \neq \emptyset; \ 1 \leq i < j \leq K\}$$

and excludes the boundary $\partial \Gamma$ by definition. An extension of the method to a variant including $\partial \Gamma$ is straightforward and not considered, neither here nor in [8]. The interfaces $\gamma_{ij}$ are numerated so that

$$\gamma = \gamma_1 \cup \ldots \cup \gamma_L.$$

Throughout we use the notation $v_i$ for the restriction of a function $v$ to the sub-domain $\Gamma_i$. Also, the jump on $\gamma$ of the function $v$ is denoted by

$$[v] : \gamma \to \mathbb{R}; \quad [v]|_{\gamma_{ij}} := v_j|_{\gamma_{ij}} - v_i|_{\gamma_{ij}} \quad (\gamma_{ij} \in \{\gamma_1, \ldots, \gamma_L\}).$$

Corresponding to the decomposition $\mathcal{T}$ there are product Sobolev spaces

$$H^s(\mathcal{T}) := H^s(\Gamma_1) \times \cdots \times H^s(\Gamma_K) \quad (s \in [0,1])$$

with usual product norm

$$\|\cdot\|^2_{H^s(\mathcal{T})} := |\cdot|^2_{H^s(\mathcal{T})} + \|\cdot\|^2_{L^2(\Gamma)}, \quad |\cdot|^2_{H^s(\mathcal{T})} := \sum_{i=1}^K |\cdot|^2_{H^s(\Gamma_i)}.$$

To define the discrete method we need some discrete spaces. On each sub-domain $\Gamma_i$ we consider quasi-uniform meshes $\mathcal{T}_{i,h_i}$ consisting of shape regular elements (parallelograms or triangles) and with mesh size $h_i$ (a precise definition is given below). We then define discrete spaces related to the decompositions $\mathcal{T}_{i,h_i}$ by

$$X_{i,h_i} := \{v \in C^0(\Gamma_i) : v|_T \in \mathbb{P}_1(T) \ \forall T \in \mathcal{T}_{i,h_i}, \ v|_{\partial\Gamma \cap \partial\Gamma_i} = 0\}.$$

Here, $\mathbb{P}_1(T)$ denotes the space of linear respectively bilinear functions on $T$ being a triangle respectively a parallelogram. We need some mesh parameters:

$$h_T := \text{diam}(T) \quad (T \in \mathcal{T}_{i,h_i}, \ i \in \{1,\ldots,K\}), \qquad h_i := \max\{h_T; \ T \in \mathcal{T}_{i,h_i}\}, \quad i = 1,\ldots,K,$$
$$h := \max\{h_i; \ i = 1,\ldots,K\}, \qquad\qquad h_{\min} := \min\{h_i; \ i = 1,\ldots,K\},$$
$$h_\mathcal{T} : \ h_\mathcal{T}|_T := h_T \ \forall T \in \mathcal{T}_{i,h_i}, \ i \in \{1,\ldots,K\}.$$

Throughout this paper we assume that $h < 1$. This is just to simplify the writing of some logarithmic terms. By $\mathcal{T}_h$ we denote the global mesh on $\Gamma$ that is defined by the sub-meshes previously defined on sub-domains. The corresponding discrete space of piecewise polynomials on $\Gamma$ is

$$X_h := X_{1,h_1} \times \cdots \times X_{K,h_K}.$$

Note that the space $X_h$ cannot be used directly to solve the discrete problem (3). $X_h$ is not a subspace of $\tilde{H}^{1/2}(\Gamma)$ since its elements are in general discontinuous across the sub-domain skeleton $\gamma$. The hypersingular operator $W$ is not well defined for discontinuous functions. Therefore, in order to utilize the non-conforming space $X_h$ for the approximation of (1) we have to rewrite its bilinear form. This amounts to applying an integration by parts formula. It involves surface differential operators which we define next. In the case of our plane surface (identified with a polygonal domain in $\mathbb{R}^2$) we simply need

$$\mathbf{curl}_\Gamma \varphi := (\partial_{x_2}\varphi, -\partial_{x_1}\varphi, 0), \qquad \text{curl}_\Gamma \boldsymbol{\varphi} := \partial_{x_1}\varphi^2 - \partial_{x_2}\varphi^1$$

for sufficiently smooth functions $\varphi$ and $\boldsymbol{\varphi} = (\varphi^1, \varphi^2, \varphi^3)$ on $\Gamma$. For an extension to Lipschitz surfaces see [3]. We also need the corresponding piecewise operators

$$\mathbf{curl}_\mathcal{T} \varphi := \sum_{i=1}^K (\mathbf{curl}_{\Gamma_i} \varphi_i)^0, \qquad \text{curl}_\mathcal{T} \boldsymbol{\varphi} := \sum_{i=1}^K (\text{curl}_{\Gamma_i} \boldsymbol{\varphi}_i)^0.$$

Here, $\mathbf{curl}_{\Gamma_i}$ and $\mathrm{curl}_{\Gamma_i}$ are the restrictions to $\Gamma_i$ of the corresponding operators on $\Gamma$, and $(\cdot)^0$ indicates extension by zero to $\Gamma$. Let $V$ be the single layer potential operator

$$V\boldsymbol{\varphi}(x) := \frac{1}{4\pi}\int_\Gamma \frac{\boldsymbol{\varphi}}{|x-y|}dS_y, \quad \boldsymbol{\varphi} \in (\tilde{H}^{-1/2}(\Gamma))^3, \quad x \in \Gamma.$$

Then there holds

$$\langle Wu, v\rangle_\Gamma = \langle \mathbf{curl}_\Gamma u, V\mathbf{curl}_\Gamma v\rangle_\Gamma \quad \text{for } u,v \in \tilde{H}^{1/2}(\Gamma),$$

cf. [26]. For discontinuous functions there is a generalization of this formula which amounts to integration by parts, see [14]. Based on this formula, a well-posed setting for discontinuous functions can be given. It uses the bilinear form

$$a(u,v) := \langle V\mathbf{curl}_\mathcal{T} u, \mathbf{curl}_\mathcal{T} v\rangle_\mathcal{T} + \langle Tu, [v]\rangle_\gamma - \langle [u], Tv\rangle_\gamma + \nu\langle [u], [v]\rangle_\gamma \qquad (4)$$

defined on $X_h \times X_h$. Here, $\nu > 0$ is a given parameter and $T$ is defined by

$$Tu := \mathbf{t}\cdot V\mathbf{curl}_\mathcal{T} u \quad \text{on} \quad \gamma$$

with $\mathbf{t}$ being a unit tangential vector on $\gamma$ which is compatible with the orientation of the jumps, cf. [8] for details. Then, the non-conforming boundary element scheme is: *Find $u_h \in X_h$ such that*

$$a(u_h, v) = \langle f, v\rangle_\Gamma \qquad \forall v \in X_h. \qquad (5)$$

According to [8, Theorem 3.1] this scheme converges almost quasi-optimally: For $\nu$ sufficiently large, and $u \in H^r(\Gamma)$ ($r \in (1, 2, 1)$) being the solution of (1), there holds

$$\sum_{j=1}^K |u - u_h|^2_{H^{1/2}(\Gamma_j)} + \nu\|[u_h]\|^2_{L^2(\gamma)} \lesssim |\log h_{\min}|^3 h^{2(r-1/2)} \|u\|^2_{H^r(\Gamma)}. \qquad (6)$$

**Remark 3.1.** Let us emphasize the following properties related to the formulation (5). They guarantee the uniqueness of the solution of the discrete method and are essential for our a posteriori error analysis. See [8, Lemmas 4.3, 4.4] and [20, Propositions 4.5, 4.6] for details.

1. Consistency: Given $f \in L^2(\Gamma)$, the solution $u$ of (1) solves (5).

2. Ellipticity of $a(\cdot, \cdot)$:

$$a(v,v) \gtrsim \|\mathbf{curl}_\mathcal{T} v\|^2_{\tilde{H}^{-1/2}(\Gamma)} + \nu\|[v]\|^2_{L^2(\gamma)} \quad \forall v \in X_h + H^r(\mathcal{T}) \quad (r > 1/2). \qquad (7)$$

3. The bilinear form defines a norm $\|v\|_a := \sqrt{a(v,v)}$ on $X_h + H^r(\mathcal{T})$ for $r > 1/2$.

# 4 A posteriori error estimation

In this section we propose and analyze an $h$-$h/2$ a posteriori error estimator.

First let us formulate the saturation assumption to be used. For a given mesh $\mathcal{T}_h$ with corresponding space $X_h$, let $\mathcal{T}_{h/2}$ denote the uniformly refined mesh with discrete space $X_{h/2}$. The corresponding discrete solutions of (5) are denoted by $u_h \in X_h$ and $u_{h/2} \in X_{h/2}$. Then the saturation assumption is that there exists $h_0 > 0$ and a constant $0 < c_{\text{sat}} < 1$ such that

$$\|u - u_{h/2}\|_a \leq c_{\text{sat}} \|u - u_h\|_a \qquad \forall h \leq h_0. \tag{8}$$

We define the following estimator terms

$$\Theta_1 := \|h_{\mathcal{T}}^{1/2} \mathbf{curl}_{\mathcal{T}} (u_h - u_{h/2})\|_{L^2(\Gamma)}, \qquad \Theta_2 := \|[u_h - u_{h/2}]\|_{L^2(\gamma)}.$$

Our main result is the following theorem.

**Theorem 4.1.** *Let $u \in \tilde{H}^{1/2}(\Gamma)$ be the solution of (1) and $u_h \in X_h$ be the solution of the discrete problem (5). Assume that (8) holds. Then, for any $\nu > 0$ and $\epsilon > 0$ there holds*

$$(1 + c_{\text{sat}})^{-1} \left( \Theta_1 + \sqrt{\nu} \Theta_2 \right) \lesssim \|u - u_h\|_a \leq C(1 - c_{\text{sat}})^{-1} h_{\min}^{-\epsilon} \Theta_1 \tag{9}$$

*with a constant $C > 0$ that depends on $\nu$ and $\epsilon$ but not on the mesh. In particular, choosing $\nu \simeq \epsilon^{-1} \simeq |\log h_{\min}|$, one has quasi-reliability of the error estimator $\Theta_1$ in the form*

$$\|u - u_h\|_a \leq C(1 - c_{\text{sat}})^{-1} |\log h_{\min}| \Theta_1 \tag{10}$$

*with a constant $C > 0$ which is independent of the mesh.*

**Remark 4.2.** Estimate (9) establishes reliability of $\Theta_1$ for $\nu > 0$ up to a slight loss of order of convergence, i.e. we only proved quasi-reliability. Selecting $\epsilon$ and $\nu$ depending on $h_{\min}$, and considering the dependence of the constant $C$ in (9) on $\epsilon$, this reveals a logarithmic perturbation in $h_{\min}$ of the reliability bound (10). We do not know whether our bound is optimal, and it is difficult to observe such perturbations in numerical experiments. On the other hand, logarithmic perturbations appear naturally in error estimates of non-conforming boundary elements, see (6).

**Remark 4.3.** Obviously, the term $\Theta_2$ in the estimate (9) can be dropped, and can also be added on the right-hand side. Therefore, both $\Theta_1$ and $\Theta_1 + \sqrt{\nu} \Theta_2$ are efficient and quasi-reliable error estimators. In our numerical experiments, the adaptive algorithm is based on $\Theta_1$ and $\Theta_2$.

## 4.1 Technical results

In this subsection we collect technical results needed for the proof of Theorem 4.1. The main ingredients are the Galerkin projector defined by the boundary element scheme and a Clement-type interpolation operator.

First let us formally define the Galerkin projector $\mathbb{G}_h$. By consistency (cf. Remark 3.1) the discrete solution $u_h \in X_h$ of (5) can be written like $u_h = \mathbb{G}_h u$ with

$$\mathbb{G}_h : H^s(\mathcal{T}) \to X_h : \quad a(\mathbb{G}_h v, w) = a(v, w) \quad \forall w \in X_h. \tag{11}$$

We start by proving an error estimate for the Galerkin approximation in the $a$-norm. This amounts to an error estimate like (6) for the non-conforming method (5) in a different norm. The proof of this bound is essentially an intermediate result from [8].

**Proposition 4.4.** *Let $s \in (1/2, 1]$ and $\psi \in H^s(\mathcal{T})$. Then, for $\nu > 0$ there holds*

$$\|\psi - \mathbb{G}_h \psi\|_a \leq D(s, \nu, h_{\min}) \inf_{v \in X_h} \|\psi - v\|_{H^s(\mathcal{T})} \tag{12}$$

*with*

$$D(s, \nu, h_{\min}) \leq C\Big(\nu^{-1/2}(s - 1/2)^{-3/2} + h_{\min}^{-(s-1/2)}(s - 1/2)^{-1} + \nu^{1/2}(s - 1/2)^{-1/2}\Big)$$

*and a positive constant $C$ which is independent of $s$, $\nu$ and the mesh.*

*Proof.* By the triangle inequality we find for any $v \in X_h$ that

$$\|\psi - \mathbb{G}_h \psi\|_a \leq \|\psi - v\|_a + \|v - \mathbb{G}_h \psi\|_a. \tag{13}$$

Discrete ellipticity (7) of $a(\cdot, \cdot)$ yields

$$\|v - \mathbb{G}_h \psi\|_a \simeq \sup_{w \in X_h \setminus \{0\}} \frac{a(v - \mathbb{G}_h \psi, w)}{\|\mathbf{curl}_\mathcal{T} w\|_{\tilde{H}^{-1/2}(\Gamma)} + \sqrt{\nu}\|[w]\|_{L^2(\gamma)}}. \tag{14}$$

According to [8, (4.18)] we can estimate

$$\frac{a(v - \mathbb{G}_h \psi, w)}{\|\mathbf{curl}_\mathcal{T} w\|_{\tilde{H}^{-1/2}(\Gamma)} + \sqrt{\nu}\|[w]\|_{L^2(\gamma)}}$$
$$\lesssim \Big((s - 1/2)^{-1} + \nu^{-1/2}(s - 1/2)^{-3/2}$$
$$\quad + h_{\min}^{-(s-1/2)}(s - 1/2)^{-1} + \nu^{1/2}(s - 1/2)^{-1/2}\Big)\|\psi - v\|_{H^s(\mathcal{T})} \tag{15}$$

for any $\nu > 0$ and $s \in (1/2, 1]$. Actually, [8] asks for $\nu \gtrsim 1$ but the estimate above holds for any $\nu > 0$. Now, by [8, (4.15), Lemma 4.1(ii)],

$$\|\mathbf{curl}_\mathcal{T}(\psi - v)\|^2_{\tilde{H}^{-1/2}(\Gamma)} \lesssim (s - 1/2)^{-2}|\psi - v|^2_{H^s(\mathcal{T})}$$

and

$$\|[\psi - v]\|^2_{L^2(\gamma)} \lesssim (s - 1/2)^{-1}\|\psi - v\|^2_{H^s(\mathcal{T})}.$$

Therefore,

$$\|\psi - v\|_a^2 = \|\mathbf{curl}_\mathcal{T}(\psi - v)\|^2_{\tilde{H}^{-1/2}(\Gamma)} + \nu\|[\psi - v]\|^2_{L^2(\gamma)}$$
$$\lesssim \Big((s - 1/2)^{-2} + \nu(s - 1/2)^{-1}\Big)\|\psi - v\|^2_{H^s(\mathcal{T})} \tag{16}$$

for any $\nu > 0$ and $s \in (1/2, 1]$. Combination of (13)–(16) proves the assertion. □

Next we introduce an interpolation operator $\mathcal{I}$ needed for the proof of Theorem 4.1. We follow the presentation of Carstensen and Bartels [4] to define operators $\mathcal{I}_i$ on $\Gamma_i$ ($i = 1, \ldots, K$). Then $\mathcal{I}$ will be simply a composition of the operators $\mathcal{I}_i$.

Let us denote by $\mathcal{N}_i$ the set of nodes of the triangulation $\mathcal{T}_{i,h_i}$, and let $\mathcal{K}_i$ be the subset of nodes which are not on the boundary $\partial\Gamma$. In the case that $\Gamma_i$ does not touch $\partial\Gamma$, $\mathcal{K}_i = \mathcal{N}_i$. For each node $z \in \mathcal{N}_i$ let $\varphi_z$ be the nodal basis function associated with this node, i.e., $\varphi_z$ is continuous, piecewise linear, $\varphi_z(x) = 0$ if $x \in \mathcal{N}_i \setminus \{z\}$ and $\varphi_z(z) = 1$. The support of $\varphi_z$ is denoted by $\omega_z$. The functions $\varphi_z$ can be modified in a straightforward way to a set of functions $\psi_z$ of the same cardinality which represents a partition of unity. In fact, in the case that $\Gamma_i$ does not touch $\partial\Gamma$ no changes are necessary. Otherwise some functions associated with nodes close to the boundary must be changed, cf. [4] for details.

Then the operators $\mathcal{I}_i$, $\mathcal{I}$ are defined by

$$\mathcal{I}_i : \begin{cases} L^1(\Gamma_i) & \to & X_{i,h_i} \\ v & \mapsto & \sum_{z \in \mathcal{K}_i} v_z \varphi_z \end{cases} \quad \text{with} \quad v_z := \frac{\int_{\omega_z} v \psi_z \, dS}{\int_{\omega_z} \varphi_z \, dS},$$

$$\mathcal{I} : L^1(\Gamma) \to X_h, \qquad \mathcal{I}|_{L^1(\Gamma_i)} := \mathcal{I}_i \quad (i = 1, \ldots, K). \tag{17}$$

**Lemma 4.5.** *Let $H_0^1(\mathcal{T})$ be the space of functions from $H^1(\mathcal{T})$ with vanishing trace on $\partial\Gamma$. For any $\delta \in (0, 1/2)$ there holds*

$$\|v - \mathcal{I}v\|_{H^s(\mathcal{T})} \lesssim \|h_\mathcal{T}^{1-s} \mathbf{curl}_\mathcal{T} v\|_{L^2(\Gamma)} \quad \forall v \in H_0^1(\mathcal{T})$$

*uniformly for all $s \in (\delta, 1/2 + \delta)$.*

*Proof.* By [4, Theorem 2.1] there holds for any sub-domain $\Gamma_i$ whose boundary has an intersection with $\partial\Gamma$ of non-zero relative measure (so that $X_{i,h_i}$ includes a homogeneous essential boundary condition)

$$\|h_\mathcal{T}^{-1}(v - \mathcal{I}v)\|_{L^2(\Gamma_i)} \lesssim |v|_{H^1(\Gamma_i)}, \quad |v - \mathcal{I}v|_{H^1(\Gamma_i)} \lesssim |v|_{H^1(\Gamma_i)} \quad \forall v \in H_D^1(\Gamma_i) \tag{18}$$

with $H_D^1(\Gamma_i)$ being the space of $H^1(\Gamma_i)$-functions with vanishing trace on $\partial\Gamma_i \cap \partial\Gamma$. Inspection of the proof of [4, Theorem 2.1] reveals that the homogeneous boundary condition on a part of the boundary is not necessary, i.e., estimates (18) hold also on sub-domains $\Gamma_i$ with $\partial\Gamma_i \cap \partial\Gamma = \emptyset$. In that case, $H_D^1(\Gamma_i) = H^1(\Gamma_i)$. Since the meshes on $\Gamma_i$ are quasi-uniform, (18) can be equivalently written as

$$\|v - \mathcal{I}v\|_{L^2(\Gamma_i)} \lesssim \|h_\mathcal{T} \mathbf{curl}_{\Gamma_i} v\|_{L^2(\Gamma_i)}, \quad \|\mathbf{curl}_{\Gamma_i}(v - \mathcal{I}v)\|_{L^2(\Gamma_i)} \lesssim \|\mathbf{curl}_{\Gamma_i} v\|_{L^2(\Gamma_i)} \quad \forall v \in H_D^1(\Gamma_i).$$

Interpolation then proves that

$$\|v - \mathcal{I}v\|_{H^s(\Gamma_i)} \lesssim \|h_\mathcal{T}^{1-s} \mathbf{curl}_{\Gamma_i} v\|_{L^2(\Gamma_i)} \quad \forall v \in H_D^1(\Gamma_i)$$

for a fixed $s \in (0, 1)$. Bounding $s$ away from 0 and 1, interpolation with the K-method gives estimates which are uniformly equivalent to the corresponding ones with Sobolev-Slobodeckij norm, cf. [2], [16, Corollary 1], so that summation proves the assertion. □

## 4.2 Proof of Theorem 4.1

First let us note that, by the saturation assumption (8) and the triangle inequality, one immediately has the two-sided estimate

$$(1 + c_{\text{sat}})^{-1}\|u_h - u_{h/2}\|_a \leq \|u - u_h\|_a \leq (1 - c_{\text{sat}})^{-1}\|u_{h/2} - u_h\|_a. \tag{19}$$

This is usually the first step in the analysis of two-level a posteriori error estimators. Note that in the case of indefinite bilinear forms the saturation parameter $c_{\text{sat}}$ enters also the lower bound, though only its boundedness is needed, cf. [1].

The lower bound in Theorem 4.1 then follows from (19) by the discrete ellipticity (7) of $a(\cdot,\cdot)$, a norm estimate from domain decomposition and an inverse estimate [16, Lemma 4],

$$(1 + c_{\text{sat}})^2 \|u - u_h\|_a^2 \geq \|u_h - u_{h/2}\|_a^2 \simeq \|\mathbf{curl}_{\mathcal{T}}(u_h - u_{h/2})\|_{\tilde{H}^{-1/2}(\Gamma)}^2 + \nu\|[u_h - u_{h/2}]\|_{L^2(\gamma)}^2$$

$$\gtrsim \sum_{i=1}^{K} \|\mathbf{curl}_{\Gamma_i}(u_h - u_{h/2})\|_{H^{-1/2}(\Gamma_i)}^2 + \nu\|[u_h - u_{h/2}]\|_{L^2(\gamma)}^2.$$

$$\gtrsim \sum_{i=1}^{K} \|h_{\mathcal{T}}^{1/2}\mathbf{curl}_{\Gamma_i}(u_h - u_{h/2})\|_{L^2(\Gamma_i)}^2 + \nu\|[u_h - u_{h/2}]\|_{L^2(\gamma)}^2.$$

$$= \|h_{\mathcal{T}}^{1/2}\mathbf{curl}_{\mathcal{T}}(u_h - u_{h/2})\|_{L^2(\Gamma)}^2 + \nu\|[u_h - u_{h/2}]\|_{L^2(\gamma)}^2.$$

Next we prove the upper bound in Theorem 4.1. For any $s \in (1/2, 1]$ let $\mathbb{G}_h : H^s(\mathcal{T}) \to X_h$ denote the Galerkin projection from (11). Note that $X_h \subset X_{h/2}$ implies $\mathbb{G}_h u_{h/2} = u_h$, hence

$$u_h - u_{h/2} = (1 - \mathbb{G}_h)(u_h - u_{h/2}).$$

Using (19) we have that

$$\|u - u_h\|_a \leq (1 - c_{\text{sat}})^{-1}\|u_h - u_{h/2}\|_a = (1 - c_{\text{sat}})^{-1}\|(1 - \mathbb{G}_h)(u_h - u_{h/2})\|_a. \tag{20}$$

By Proposition 4.4 with $v = \mathcal{I}(u_h - u_{h/2})$ (recall the definition of $\mathcal{I}$ from (17)) and making use of Lemma 4.5 we obtain with $s = 1/2 + \epsilon$ (for $\epsilon \in (0, 1/4]$)

$$\|(1 - \mathbb{G}_h)(u_h - u_{h/2})\|_a^2 \lesssim D(s, \nu, h_{\min})^2 \|u_h - u_{h/2} - \mathcal{I}(u_h - u_{h/2})\|_{H^s(\mathcal{T})}^2$$

$$\lesssim D(s, \nu, h_{\min})^2 \|h_{\mathcal{T}}^{1-s}\mathbf{curl}_{\mathcal{T}}(u_h - u_{h/2})\|_{L^2(\Gamma)}^2$$

$$\lesssim D(1/2 + \epsilon, \nu, h_{\min})^2 h_{\min}^{-2\epsilon} \|h_{\mathcal{T}}^{1/2}\mathbf{curl}_{\mathcal{T}}(u_h - u_{h/2})\|_{L^2(\Gamma)}^2.$$

Combining this estimate with (20) (and writing $\epsilon$ instead of $2\epsilon$) proves the upper bound in Theorem 4.1, with constant $C > 0$ depending on $\nu$ and $\epsilon$. Now choosing $\epsilon = |\log h_{\min}|^{-1}$, $\nu \simeq \epsilon^{-1}$ and recalling that

$$D(1/2 + \epsilon, \nu, h_{\min}) \lesssim \nu^{-1/2}\epsilon^{-3/2} + h_{\min}^{-\epsilon}\epsilon^{-1} + \nu^{1/2}\epsilon^{-1/2},$$

cf. Proposition 4.4, this gives the specified reliability.

## 4.3 Error estimation in a standard norm

Theorem 4.1 states efficiency and quasi-reliability of the error estimator in the $a$-norm. This norm is the natural one stemming from the non-conforming (discrete) variational formulation. In this section we briefly indicate what can be proved when replacing the $a$-norm with a standard broken Sobolev norm defined by

$$\|v\|^2_{H^s_\nu(\mathcal{T})} := \sum_{i=1}^K |v|^2_{H^s(\Gamma_i)} + \nu\|[v]\|^2_{L^2(\gamma)}.$$

This norm is well defined for functions $v \in H^r(\mathcal{T})$ when $r \geq s$ and $r > 1/2$. In particular, when $s = 1/2$, it is well defined for elements of $X_h$ and the solution $u$ of (1).

**Corollary 4.6.** *Let $u \in \tilde{H}^{1/2}(\Gamma)$ be the solution of (1) and $u_h \in X_h$ be the solution of the discrete problem (5), and assume that (8) holds.*
*For any $\epsilon > 0$ (sufficiently small) and any $\nu > 0$ there holds efficiency in the sense that*

$$(1 + c_{\text{sat}})^{-1}\Big(\Theta_1 + \sqrt{\nu}\Theta_2\Big) \lesssim \epsilon^{-1}|u - u_h|_{H^{1/2+\epsilon}(\mathcal{T})} + \sqrt{\nu}\|[u - u_h]\|_{L^2(\gamma)} \leq \epsilon^{-1}\|u - u_h\|_{H^{1/2+\epsilon}_\nu(\mathcal{T})}.$$

*Furthermore, we have the reliability bound*

$$\|u - u_h\|_{H^{1/2}_\nu(\mathcal{T})} \leq C(1 - c_{\text{sat}})^{-1} h_{\min}^{-\epsilon} \Theta_1$$

*for any $\nu > 0$ and $\epsilon > 0$ with a constant $C > 0$ that depends on $\nu$ and $\epsilon$ but not on the mesh. Choosing $\nu \simeq \epsilon^{-1} \simeq |\log h_{\min}|$, one has quasi-reliability of the error estimator $\Theta_1$ in the form*

$$\|u - u_h\|_{H^{1/2}_\nu(\mathcal{T})} \leq C(1 - c_{\text{sat}})^{-1} |\log h_{\min}| \Theta_1$$

*with a constant $C > 0$ which is independent of the mesh.*

*Proof.* To show the efficiency bound one uses the continuity of $\mathbf{curl}_{\Gamma_i}$ and a quotient space argument [18]:

$$\|\mathbf{curl}_\mathcal{T}(u - u_h)\|^2_{\tilde{H}^{-1/2}(\Gamma)} \lesssim \epsilon^{-2} \sum_{i=1}^K \|\mathbf{curl}_{\Gamma_i}(u - u_h)\|^2_{H^{-1/2+\epsilon}(\Gamma_i)}$$

$$\lesssim \epsilon^{-2} \sum_{i=1}^K |u - u_h|^2_{H^{1/2+\epsilon}(\Gamma_i)}$$

for any $\epsilon > 0$, cf. [8, (4.15)] for details. Then, efficiency of the estimator $\Theta_1 + \sqrt{\nu}\Theta_2$ with respect to the $a$-norm by Theorem 4.1 and the identity

$$a(u - u_h, u - u_h) = \|\mathbf{curl}_\mathcal{T}(u - u_h)\|^2_{\tilde{H}^{-1/2}(\Gamma)} + \nu\|[u - u_h]\|^2_{L^2(\gamma)}$$

prove the efficiency estimate of the corollary.

Now we prove quasi-reliability. By the $\tilde H^{-1/2}(\Gamma)$-ellipticity of $V$ and a standard norm estimate we have that

$$\langle V\mathbf{curl}_{\mathcal{T}}v, \mathbf{curl}_{\mathcal{T}}v\rangle_{\mathcal{T}} \gtrsim \|\mathbf{curl}_{\mathcal{T}}v\|^2_{\tilde H^{-1/2}(\Gamma)} \gtrsim \sum_{i=1}^{K} \|\mathbf{curl}_{\Gamma_i}v\|^2_{H^{-1/2}(\Gamma_i)}$$

for any $v$ with $v|_{\Gamma_i} \in H^r(\Gamma_i)$ for $r > 1/2$. By [8, Lemma 4.1],

$$\|\mathbf{curl}_{\Gamma_i}v\|_{H^{-1/2}(\Gamma_i)} \gtrsim |v|_{H^{1/2}(\Gamma_i)}$$

for any $v \in H^{1/2}(\Gamma_i)$. Combination of both estimates proves that

$$a(v,v) \gtrsim \sum_{i=1}^{K} |v|^2_{H^{1/2}(\Gamma_i)} + \nu\|[v]\|^2_{L^2(\gamma)} = \|v\|^2_{H^{1/2}_\nu(\mathcal{T})}$$

for any $v$ with $v|_{\Gamma_i} \in H^r(\Gamma_i)$ and $r > 1/2$. The quasi-reliability of $\Theta_1$ with respect to the $a$-norm by Theorem 4.1 proves the quasi-reliability of this estimator with respect to the $\|\cdot\|_{H^{1/2}_\nu(\mathcal{T})}$-norm. □

### 4.4 Adaptive refinement strategy

Based on the a posteriori error estimate (9), we define local error indicators

$$\theta_T^2 := \|h_T^{1/2}\mathbf{curl}_{\mathcal{T}}(u_h - u_{h/2})\|^2_{L^2(T)} + \nu\|[u_h - u_{h/2}]\|^2_{L^2(\gamma \cap \partial T)} \tag{21}$$

for all elements $T \in \mathcal{T}_{i,h_i}$ ($i = 1,\ldots,K$) and use a standard adaptive refinement algorithm as follows.

**Algorithm.** Let $\mathcal{T}_h^{(0)}$ be an initial mesh and set $j = 0$. For a given tolerance tol $> 0$ and a refinement parameter $\delta \in (0,1]$ perform the following steps.

1. Refine $\mathcal{T}_h^{(j)}$ uniformly to obtain $\mathcal{T}_{h/2}^{(j)}$.

2. Compute $u_h$ and $u_{h/2}$.

3. Compute the local indicators $\theta_T$ for all elements $T$.

4. Compute the error estimator

$$\Theta^2 = \sum_{i=1}^{K} \sum_{T \in \mathcal{T}_{i,h_i}} \theta_T^2. \tag{22}$$

5. Stop if $\Theta \le$ tol, otherwise continue.

6. Determine a minimal set $\mathcal{R}$ of elements $T$ with largest indicators such that

$$\sum_{T \in \mathcal{R}} \theta_T^2 \ge \delta^2 \Theta^2.$$

7. Construct a new mesh $\mathcal{T}_h^{j+1}$ by refining the elements of the set $\mathcal{R}$ and by introducing further edges subject to a minimum angle condition such that meshes on sub-domains are conforming.

8. Update $j := j+1$ and go to 1.

# 5 Numerical experiments

In this section we report on some numerical experiments to study efficiency and reliability of the a posteriori error estimator $\Theta$ defined by (22), (21) and to check the performance of the corresponding adaptive refinement strategy. To be precise, for ease of programming, the jump terms on edges of the interface between sub-domains are taking into account twice.
We consider the model problem (1) with $f = 1$ on $\Gamma = (-1/2, 1/2)^2$ and use triangular meshes $\mathcal{T}_h$. Since the exact solution $u$ of (1) is unknown, the error in standard broken Sobolev norm

$$\|u - u_h\|^2_{H_\nu^{1/2}(\mathcal{T})} = |u - u_h|^2_{H^{1/2}(\mathcal{T})} + \nu \|[u_h]\|^2_{L^2(\gamma)}$$

cannot be computed directly, except for $\|[u_h]\|_{L^2(\gamma)}$ which is straightforward to implement. As in previous papers on non-conforming approximations of hypersingular operators, see [8, 20], we take

$$(\text{``total error''})^2 := \big|\|u\|^2_{\text{ex}} - \langle f, u_h \rangle_\Gamma\big| + \nu\|[u_h]\|^2_{L^2(\gamma)} \tag{23}$$

as a computable and reasonable measure for an upper bound of $\|u - u_h\|_{H_\nu^{1/2}(\mathcal{T})}$. Here, $\|u\|_{\text{ex}}$ is an extrapolated value of $\|u\|_{\tilde{H}^{1/2}(\Gamma)}$, cf. [12].
Throughout we consider the decomposition of $\Gamma$ into four sub-domains and initial mesh as shown in Figure 1. Figure 2 shows the "total error" (23) and "total estim" $:= \Theta$ for uniform and adaptive refinements ($\delta = 0.5$, $\nu = 100$). The curves $N^{-1/4}$ and $N^{-1/2}$ with $N = \dim X_h$ are also given. They correspond to $h^{1/2}$ and $h$, respectively, for uniform meshes. As can be seen, the error and estimator curves are parallel to $h^{1/2}$ for uniform meshes. This is the expected convergence order, cf. [27]. In the adaptive case both curves appear to approach $N^{-1/2}$ for higher dimension, as expected. In Figures 3 and 4 the individual terms

$$\text{error1} := \big|\|u\|^2_{ex} - \langle f, u_h \rangle_\Gamma\big|^{1/2}, \qquad \text{error2} := \sqrt{\nu}\|[u_h]\|_{L^2(\gamma)},$$
$$\text{estim1} := \Theta_1 = \|h_\mathcal{T}^{1/2}\mathbf{curl}_\mathcal{T}(u_h - u_{h/2})\|_{L^2(\Gamma)}, \qquad \text{estim2} := \Theta_2 = \sqrt{\nu}\|[u_h - u_{h/2}]\|_{L^2(\gamma)}$$

are shown, for uniform and adaptively refined meshes, respectively. Except for "error1"/"estim1" in the adaptive case, the corresponding error and estimator curves are parallel, thus confirming reliability and efficiency. Our interpretation of the exception (where "error1" has slightly faster convergence than "estim1") is that the error calculation via (23) has a pre-asymptotic behavior which does not show the correct slope of the error in energy norm. This has been observed already in [8]. Indeed, the estimator curves show the correct behavior in all our experiments.

In Fig. 5 two samples from the sequence of adaptively refined meshes are plotted. In our case, $\Gamma$ is an open surface where the solution has strong singularities along the edges. The error indicators detect these singularities and the adaptive algorithm refines accordingly.

To see the influence of $\nu$, Figure 6 shows the errors of uniform and adaptive refinements for $\nu = 10$. Again, the uniform version converges like $h^{1/2} \simeq N^{-1/4}$ and adaptivity yields convergence of the order $N^{-1/2}$. Though, in this case of smaller $\nu$, the discrepancy between error and estimator is larger, see Figures 7 and 8 where the individual terms error1, error2, $\Theta_1$ and $\Theta_2$ are shown. Again, we attribute this to the limitation of the error "calculation" via (23).

To resume, the estimator shows the expected behavior for uniform as well as for adaptively refined meshes, for larger and smaller $\nu$. Also, meshes are adaptively refined in the expected way to resolve the edge singularities of the exact solution. Here, discontinuities of the approximation at interfaces do not have a negative influence. Indeed, mesh refinement seems not to take into account interfaces, as it should be.

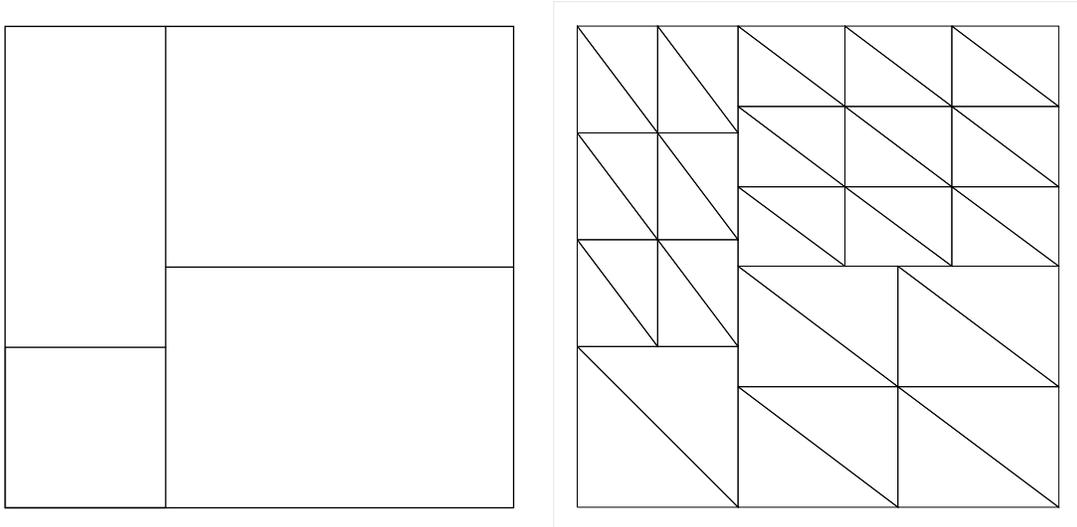

Figure 1: Domain decomposition $\mathcal{T}$ of $\Gamma$ and initial mesh.

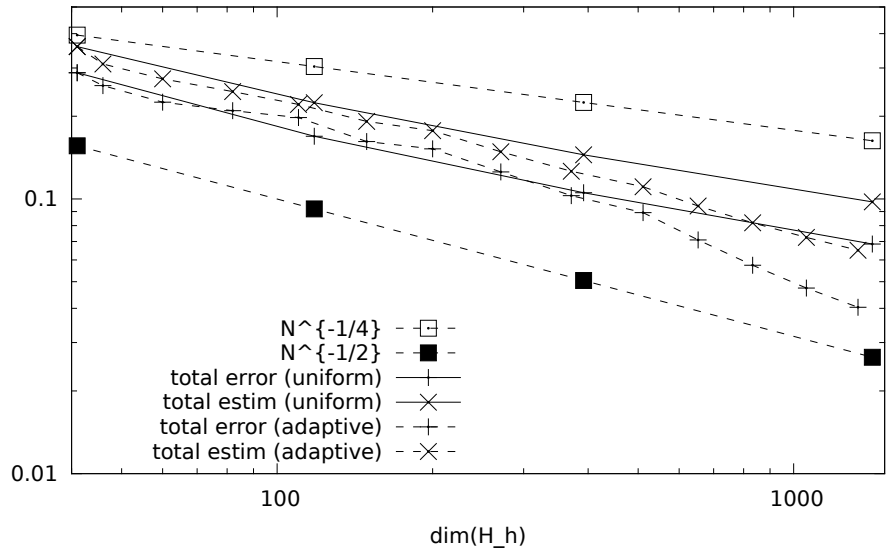

Figure 2: Total error and estimator for uniform and adaptive refinements ($\nu = 100$).

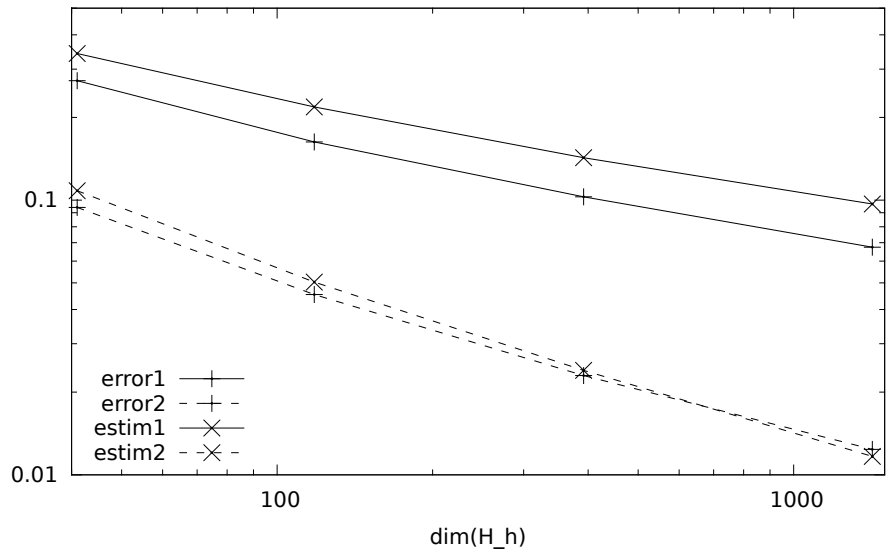

Figure 3: Individual errors and estimators for uniform refinement ($\nu = 100$).

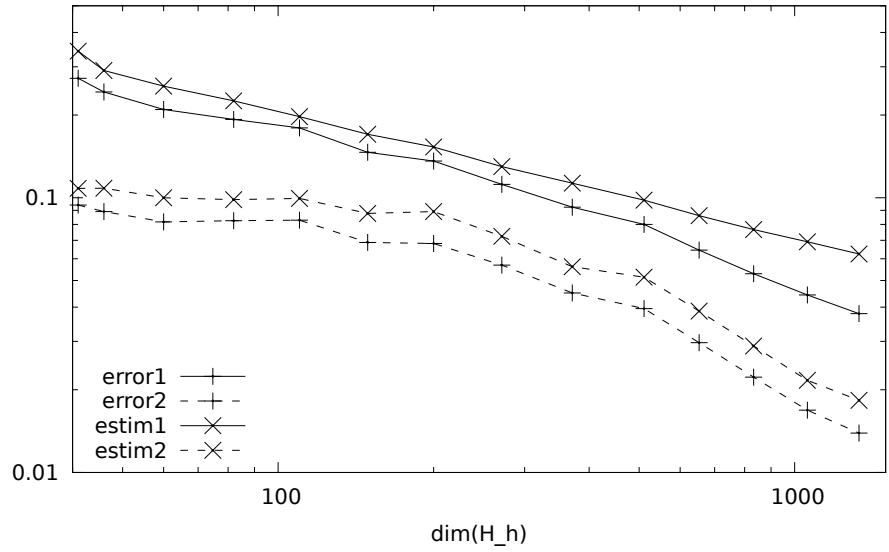

Figure 4: Individual errors and estimators for adaptive refinement ($\nu = 100$).

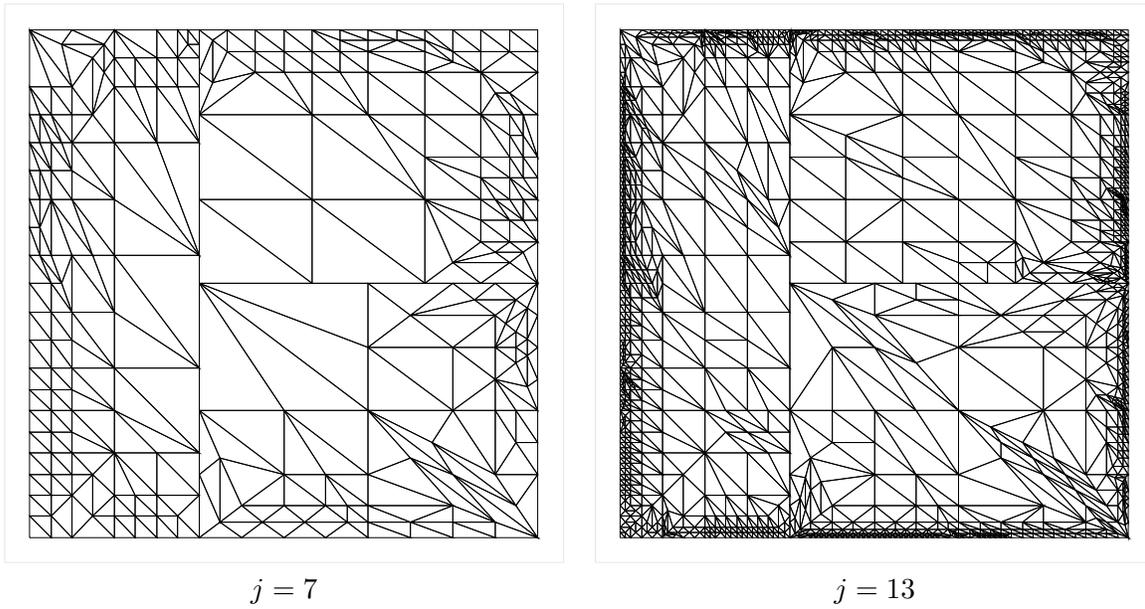

Figure 5: Adaptively refined meshes (steps 7 and 13) with parameters $\nu = 100$, $\delta = 0.5$.

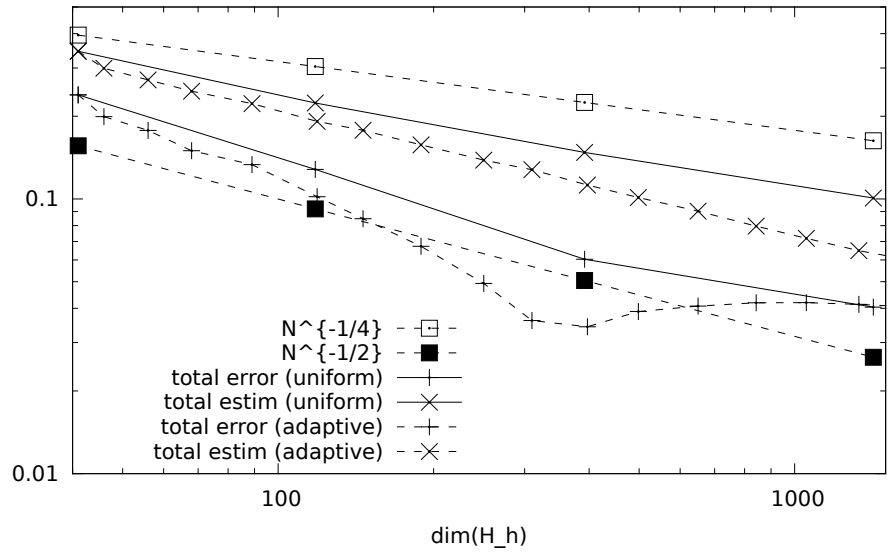

Figure 6: Total error and estimator for uniform and adaptive refinements ($\nu = 10$).

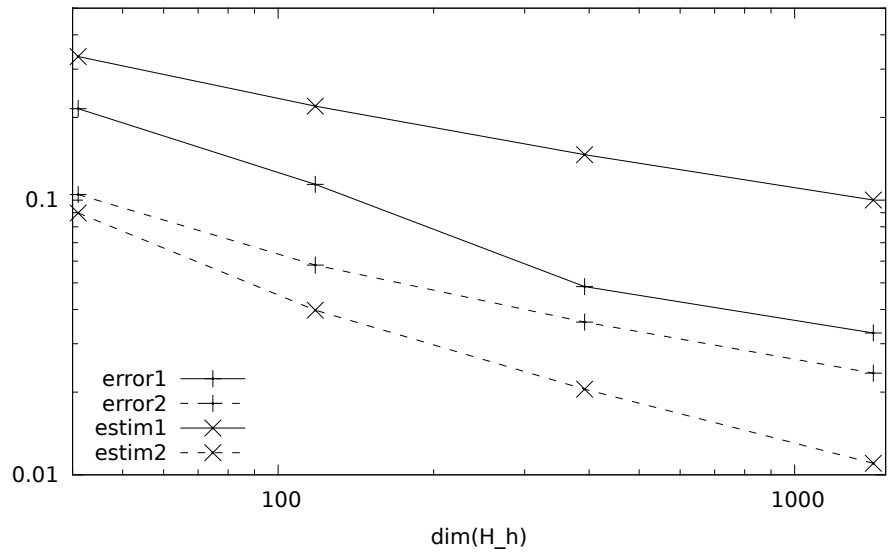

Figure 7: Individual errors and estimators for uniform refinement ($\nu = 10$).

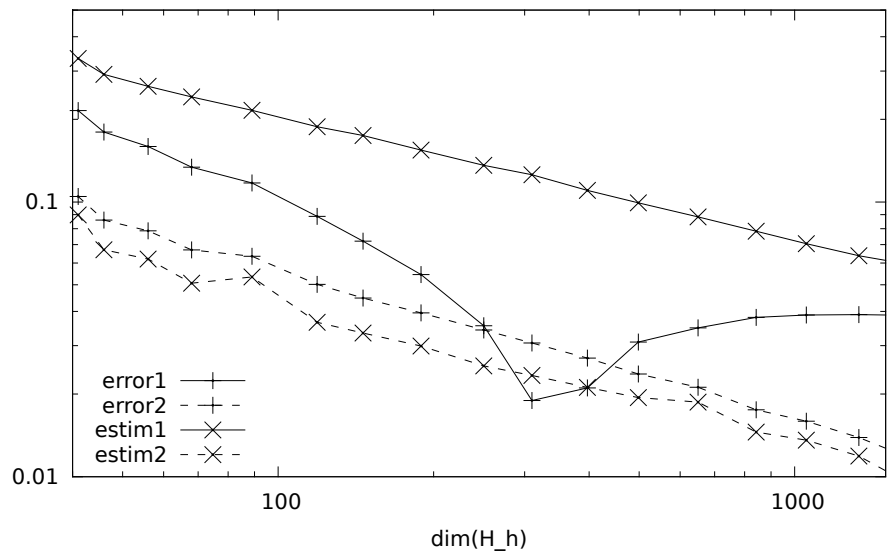

Figure 8: Individual errors and estimators for adaptive refinement ($\nu = 10$).